\documentclass[11pt]{article}

\def\sign{\mathop{\rm sign}\nolimits}%

\newtheorem{theorem}{Theorem}

\newenvironment{proof}[1]{{\bf #1:}}{\hspace*{\fill}$\Box$\newline}

\newcommand{\ind}{{\perp\hspace*{-2.1mm}\perp}}

\usepackage{apacite,latexsym,graphics,epsfig,subfigure,color,setspace,lineno}
\usepackage[left=3.25cm,right=3.25cm]{geometry}

 \title{Nonparametric testing of conditional independence by means of the partial copula
}%
 \author{Wicher P.\ Bergsma\\London School of Economics and Political Science}

\begin{document}

\def\prob{\mathop{\rm P}\nolimits}%
\newcommand{\nind}{{\not\perp\hspace*{-2.1mm}\perp}}



\maketitle

 \begin{abstract}
We propose a new method to test conditional independence of two
real random variables $Y$ and $Z$ conditionally on an arbitrary third random variable $X$.
The partial copula is introduced, defined as the joint distribution of $U=F_{Y|X}(Y|X)$ and $V=F_{Z|X}(Z|X)$.
We call this transformation of $(Y,Z)$ into $(U,V)$ the partial copula transform.
It is easy to show that if $Y$ and $Z$ are continuous for any given value of $X$, then $Y\ind Z|X$ implies $U\ind V$. Conditional independence can then be tested by (i) applying the partial copula transform to the data points and (ii) applying a test of ordinary independence to the transformed data. In practice, $F_{Y|X}$ and $F_{Z|X}$ will need to be estimated, which can be done by, e.g., standard kernel methods. We show that under easily satisfied conditions, and for a very large class of test statistics for independence which includes the covariance, Kendall's tau, and Hoeffding's test statistic, the effect of this estimation vanishes asymptotically.
Thus, for large samples, the estimation can be ignored and we have a simple method which can be used to apply a wide range of tests of independence, including ones with consistency for arbitrary alternatives, to test for conditional independence.
A simulation study indicates good small sample performance.
Advantages of the partial copula approach compared to competitors seem to be simplicity and generality.
\end{abstract}%


\newcommand{\X}{\mbox{\rm X}}

\section{Introduction}  \label{intro}

Random variables $Y$ and $Z$ are conditionally independent given $X$ if, knowing the value of $X$, information about $Y$ does not provide information about $Z$.
Following \citeA{dawid79}, this hypothesis is denoted as $Y\ind Z|X$.
Conditional independence models are the building blocks of graphical models, which can be used for causal modeling \cite{wc96,lauritzen96,pearl00}. The literature on the topic is usually restricted to the normal and categorical cases, but recently nonparametric testing has also received a fair amount of attention \cite{sw07,sw08,song09,huang10,brt10}. These papers also contain useful further references and applications. At the end of this section a brief literature review is presented as well. This paper presents the partial copula approach for the nonparametric testing of the conditional independence hypothesis for real random variables $Y$ and $Z$ controlling for an arbitrary random variable $X$ based on $n$ independent and identically distributed (iid) replications of $(X,Y,Z)$.

In general, compared to the testing of ordinary independence, the present problem seems a difficult one (\citeNP{bergsma04}). In particular, if, without further assumptions, for a given value of $X$ only one $(Y,Z)$ observation is available, this does not tell us anything at all about the conditional association for that particular value of $X$. Alternatively, it can be said that the difficulty lies in the possible presence of marginal dependencies of $Y$ on $X$ and of $Z$ on $X$: if both $Y\ind X$ and $Z\ind X$, then $Y\ind Z|X$ implies $Y\ind Z$. This is easily shown as follows: assuming the conditions hold,
\begin{eqnarray}
   \lefteqn{P(X\in A,Y\in B) = P(X\in A,Y\in B|X=x)}\nonumber\\
   && = P(X\in A|X=x)P(Y\in B|X=x) = P(X\in A)P(Y\in B)\label{eq1}
\end{eqnarray}
for measurable subsets $A$ and $B$ of the sample spaces of $Y$ and $Z$, an $x$ in the sample space of $X$.
Thus, in the absence of the marginal associations, conditional independence for $(X,Y,Z)$ can be tested simply by performing a test of ordinary independence for $(Y,Z)$.
In practice, of course, marginal dependencies are typically present, and conditional independence might be tested by removing (some of) the marginal dependency.
Below, we describe two methods which use this idea.

Firstly, the idea of removing some of the marginal association is used to construct the well-known partial correlation coefficient.
The assumption is made that  $Y=g(X)+U'$ and $Z=h(X)+V'$ where the error terms $U'$ and $V'$ have zero expectation conditionally on $X$. Now clearly the correlations between $X$ and $U'$ and between $X$ and $V'$ are zero, i.e., replacing $Y$ by $U'$ and $Z$ by $V'$ removes some of the marginal association. The partial correlation coefficient is the correlation between $U'$ and $V'$, which is zero under conditional independence. Thus, the partial correlation coefficient can be used for a nonparametric test of conditional independence.

Secondly, we present the approach of this paper, which has two key distinguishing features compared to the partial correlation approach: (i) the use of conditional ranks and (ii) complete removal of the marginal dependence.
It is based on the new concept of the partial copula, which, like the partial correlation, is easily defined.
Let $X$ be a random variable on an arbitrary measurable space and let $Y$ and $Z$ be real-valued random variables with
conditional distribution functions $F_{Y|X}$ and $F_{Z|X}$. With
\begin{eqnarray}
   U = F_{Y|X}(Y|X)
   \hspace{6mm}\mbox{and}\hspace{6mm}
   V = F_{Z|X}(Z|X)\label{pct}
\end{eqnarray}
the partial copula pertaining to $Y$ and $Z$ given $X$ is defined as the joint distribution function of $(U,V)$ (cf.\ the ordinary copula pertaining to $(Y,Z)$ is the joint distribution function of $(Y,Z)$).
We call transformation~(\ref{pct}) the partial copula transform.
If both $Y$ and $Z$ have continuous conditional distribution functions for every given value of $X$, then both $U|X=x$ and $V|X=x$ are uniformly distributed for all $x$, hence $U\ind X$ and $V\ind X$. Therefore, by~(\ref{eq1}), $Y\ind Z|X$ implies $U\ind V$.
If we now apply the partial copula transform to iid data $(X_1,Y_1,Z_1),\ldots,(X_n,Y_n,Z_n)$, i.e., setting
\begin{eqnarray}
   U_i=F_{Y|X}(Y_i|X_i)
   \hspace{6mm}\mbox{and}\hspace{6mm}
   V_i=F_{Z|X}(Z_i|X_i),
   \label{trans1}
\end{eqnarray}
any test of ordinary independence for the $(U_i,V_i)$ is a test for conditional independence for the $(Y_i,Z_i)$ given the $X_i$.

In practice, $F_{Y|X}$ and $F_{Z|X}$ are typically unknown and need to be estimated. In Section~\ref{sec-ex} it is shown that for a wide variety of estimators, the asymptotic distribution of a broad class of test statistics in the form of a covariance between two functions on the marginals is unaffected by the estimation. Thus, a comprehensive range of tests of ordinary independence lead to an asymptotic test of conditional independence.
The procedure is illustrated using a real data set and several tests of independence, including ones based on the covariance and Kendall's tau, and some tests with consistency against any alternative such as Hoeffding's test. Questions for further research are addressed in Section~\ref{sec-rem}.

Note that, since $U\ind V$ does not imply $Y\ind Z|X$, a test of the hypothesis $U\ind V$ cannot be consistent for all alternatives to the hypothesis $Y\ind Z|X$: even if $U\ind V$ there may still be so-called three-variable interaction present. Fortunately, the testing for three-variable interaction seems in general a lot easier than the testing for conditional independence \cite{bkr61}, and it is reasonable to use a separate procedure for that.
Using the partial copula transform, we can expect most power against alternatives with little three-variable interaction, in particular, alternatives for which the conditional copula, i.e., the joint distribution of $(F_{Y|X}(Y|x),F_{Z|X}(Z|x))$, does not vary much with~$x$.

The following are some recent approaches with the same aim as the present paper. \citeA{sw07,sw08} compare the conditional densities $f_{Y|XZ}$ and $f_{Y|X}$ using characteristic functions and Hellinger distances, while \citeA{brt10} compare the conditional distribution functions $F_{Y|XZ}$ and $F_{Y|X}$. Kernel methods are used to estimate the densities and distribution functions.
The methods of \citeA{song09} and \citeA{huang10} are closer to the present paper's, in that they do not require estimation of the distribution of $Y$ given both $X$ and $Z$, but only given $X$. Huang uses R{\'e}nyi type maximal correlations and differs from our approach in that it is not based on conditional ranks. Song's also uses the transformation~(\ref{trans1}), which he calls the Rosenblatt transform. Song and the present author found this independently, the present paper having appeared in different form as a technical report \cite{bergsma04}.
\citeA{lg97} considered conditional independence tests based on Cram{\'e}r von Mises and Kolmogorov-Smirnov criteria.

A brief selection of older approaches is as follows.
Kendall \citeyear{kendall42} introduced a partial version of Kendall's tau. However it is not clear how useful this coefficient can be in practice, as it is not
necessarily zero under conditional independence unless certain restrictive conditions are met \cite{korn84}.
\citeA{goodman59} and \citeA{gripenberg92}) proposed a method based on another partial version of Kendall's tau, using the number of local concordant and discordant pairs of observations.
The easiest case for testing conditional independence is if $X$ is categorical, with sufficiently many
observations per category: for each category, a test of independence
can be done, and these tests can be combined in various ways. If all
three variables are categorical, log-linear techniques can be used
(cf.\ \citeNP{lauritzen96}, \citeNP{agresti02}).

\section{A practical approach to conditional independence testing}\label{sec-ex}

If $F_{Y|X}$ and $F_{Z|X}$ are known, the partial copula transform~(\ref{trans1}) can be applied and an arbitrary test of ordinary independence can be applied to the $(U_i,V_i)$. A permutation test would then be valid even for small samples. In Section~\ref{sec-ind} a general class of test statistics for independence testing is described. In practice, however, $F_{Y|X}$ and $F_{Z|X}$ are usually unknown and need to be estimated; in Section~\ref{sec-ht}, we show that, under fairly unrestrictive conditions, this estimation does not affect the asymptotic behaviour of the test statistics given in Section~\ref{sec-ind}.
In Section~\ref{sec-digox}, we illustrate the procedure on a real-data example using several tests of independence and provide some graphical displays.

\subsection{A class of test statistics for ordinary independence}\label{sec-ind}

Many measures of association are of the form
\begin{eqnarray}
   \theta(Y,Z) = Es(Y_1,\ldots,Y_r)t(Z_1,\ldots,Z_r)\label{theta}
\end{eqnarray}
for some $r\ge 1$ and where $(Y_i,Z_i)$ are independent replications of $(Y,Z)$. Square integrability of $s$ and $t$ is a sufficient condition for the expectation to exist.
If  $s$ and $t$ have zero means, then $Y\ind Z$ implies $\theta=0$. The sample ($V$-statistic) estimator which can be used as a test statistic for independence is
\begin{eqnarray}
   \hat\theta(Y,Z) = \frac{1}{n^r}\sum_{i_1=1}^n\ldots\sum_{i_r=1}^n s(Y_{i_1},\ldots,Y_{i_r})t(Z_{i_1},\ldots,Z_{i_r})\label{vstat}
\end{eqnarray}
In the absence of distributional assumptions the best way to compute a $p$-value based on $\hat\theta$ is usually the permutation test, but alternatively asymptotic theory is well-developed in certain cases \cite{hoeffding48,rw79,serfling80}.

The ordinary covariance and Kendall's tau are well-known association measures which can be written in the form~(\ref{theta}).
Spearman's rho can of course be written in that form as well, but we need not consider it because for rank data it is proportional to the covariance.
Some coefficients of form~(\ref{theta}) which are nonnegative and zero iff independence holds are the following.
The first is Hoeffding's $\Delta$, which can be defined as
\begin{eqnarray*}
   \Delta(Y,Z) = \frac14 E\phi((Y_1,Y_2,Y_3)\phi((Y_1,Y_4,Y_5)\phi((Z_1,Z_2,Z_3)\phi((Z_1,Z_4,Z_5)
\end{eqnarray*}
where $\phi(z_1,z_2,z_3)=I(z_1\ge z_2)-I(z_1\ge z_3)$.
\citeA{bergsma06} introduced
\begin{eqnarray*}
    \kappa(Y,Z)
    =
    \frac14Ea(Y_1,Y_2,Y_3,Y_4)a(Z_1,Z_2,Z_3,Z_4)
\end{eqnarray*}
and
\begin{eqnarray*}
    \tau^*(Y,Z)
    =
    E\sign[a(Y_1,Y_2,Y_3,Y_4)]\sign[a(Z_1,Z_2,Z_3,Z_4)]
\end{eqnarray*}
where
\begin{eqnarray*}
    a(z_1,z_2,z_3,z_4) = |z_1-z_2|+|z_3-z_4|-|z_1-z_3|-|z_2-z_4|
\end{eqnarray*}
The latter is an extension of Kendall's tau introduced by \citeA{bd10}.
Note that $\Delta$, $\kappa$ and $\tau^*$ lead to tests consistent for all alternatives to independence.


\subsection{Using estimated conditional distributions}\label{sec-ht}

In practice, $F_{Y|X}$ and $F_{Z|X}$ are typically unknown and need to be estimated. With estimates $\hat F_{Y|X}$ and $\hat F_{Z|X}$, the partial copula transform yields
\begin{eqnarray}
   \hat U_i = \hat F_{Y|X}(Y_i|X_i)
   \hspace{6mm}\mbox{and}\hspace{6mm}
   \hat V_i = \hat F_{Z|X}(Z_i|X_i)\label{newobs}
\end{eqnarray}
and substitution into~(\ref{vstat}) gives
\begin{eqnarray*}
   \hat\theta(\hat U,\hat V) = \frac{1}{n^r}\sum_{i_1=1}^n\ldots\sum_{i_r=1}^n s(\hat U_{i_1},\ldots,\hat U_{i_r})t(\hat V_{i_1},\ldots,\hat V_{i_r})
\end{eqnarray*}
Theorem~\ref{th1} states that under the following easily satisfied conditions, the fact that we are using estimates $(\hat U_i,\hat V_i)$ rather than true values $(U_i,V_i)$ can be ignored.
\newcounter{regu}
\begin{list}
   {A\arabic{regu}:}{\usecounter{regu}}
\item $s$ and $t$ are continuous almost everywhere with respect to Lebesgue measure on $[0,1]^r$,\label{a1}

\item $n^\alpha[\hat\theta(U,V)-\theta(U,V)]=O_p(1)$ for some positive $\alpha$,\label{a2}

\item for almost all $(x,y,z)$ and some positive $\beta_1$ and $\beta_2$,\label{a3}
\begin{eqnarray*}
   n^{\beta_1}[\hat F_{Y|X}(y|x)- F_{Y|X}(y|x)] = O_p(1)\\
   n^{\beta_2}[\hat F_{Z|X}(z|x)- F_{Z|X}(z|x)] = O_p(1)
\end{eqnarray*}
with uniform convergence on all compact sets.
\end{list}
\begin{theorem}\label{th1}
Assume {\rm A1}, {\rm A2}, and {\rm A3} hold.  Then
\begin{eqnarray*}
   n^{\alpha}[\hat\theta(\hat U,\hat V) - \theta(U,V)] = n^{\alpha}[\hat\theta(U,V) - \theta(U,V)] + O_p\left(n^{-\min(\beta_1,\beta_2)}\right)
\end{eqnarray*}
\end{theorem}
\proof{Proof}
Continuity of $s$ and $t$ (a.e.) implies continuity of $s\times t$ (a.e.). Hence,
\begin{eqnarray*}
   \lefteqn{s(\hat U_{i_1},\ldots,\hat U_{i_r})t(\hat V_{i_1},\ldots,\hat V_{i_r}) - \theta(U,V)
   =}\\
   &&
   \Big(s(U_{i_1},\ldots,U_{i_r})t(V_{i_1},\ldots,V_{i_r}) - \theta(U,V)\Big)\left(1 + O_p\left(n^{-\min(\beta_1,\beta_2)}\right)\right)
\end{eqnarray*}
with probability 1.
From this, the uniform convergence, and A2,
\begin{eqnarray*}
   \lefteqn{n^{\alpha}[\hat\theta(\hat U,\hat V) - \theta(U,V)]}\\
   &=& n^{\alpha}\left(\frac{1}{n^r}\sum_{i_1=1}^n\ldots\sum_{i_r=1}^n s(\hat U_{i_1},\ldots,\hat U_{i_r})t(\hat V_{i_1},\ldots,\hat V_{i_r}) - \theta(U,V)\right)\\
   &=& n^{\alpha}\frac{1}{n^r}\sum_{i_1=1}^n\ldots\sum_{i_r=1}^n\Big[s(U_{i_1},\ldots,U_{i_r})t(V_{i_1},\ldots,V_{i_r}) - \theta(U,V)\Big]\left[1 + O_p\left(n^{-\min(\beta_1,\beta_2)}\right)\right]\\
   &=& n^{\alpha}[\hat\theta(U,V) - \theta(U,V)]\left[1 + O_p\left(n^{-\min(\beta_1,\beta_2)}\right)\right]\\
   &=& n^{\alpha}[\hat\theta(U,V) - \theta(U,V)] + O_p\left(n^{-\min(\beta_1,\beta_2)}\right)
\end{eqnarray*}
with probability 1.
\endproof
Note that the proof is immediately adapted to $U$-statistic estimators. Assumptions A1 and A2 are satisfied for the measures of association mentioned in Section~\ref{sec-ind}, i.e., the covariance, Kendall's tau, and $\Delta$, $\kappa$ and $\tau^*$. For Euclidean $X$, Nadaraya-Watson estimators are simple estimators of the marginal distribution functions satisfying A3 (under some additional restrictions).
Alternatively, \citeA{hwy99} introduced local linear and local logistic estimators, which improve the Nadaraya-Watson ones in certain respects.
However, local linear estimators may not lead to a distribution function, and \citeauthor{hwy99} developed adjusted estimators which solve this problem.

\subsection{Application with graphical illustration}\label{sec-digox}

\begin{figure}[tbp]
\centering
 \subfigure[Untransformed marginal distributions with association present]
  {\includegraphics[width=.85\linewidth,clip]{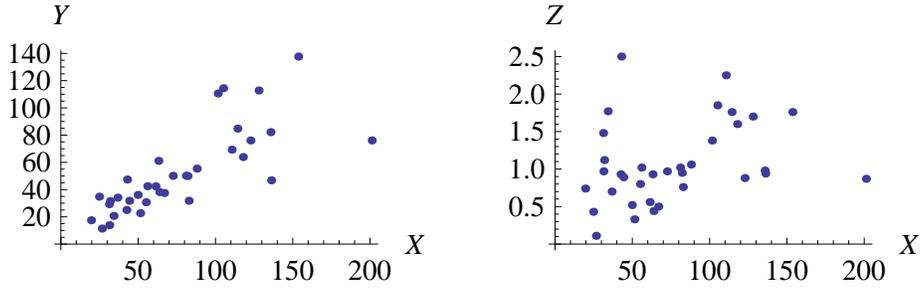}
}
 \hspace{8mm}
 \subfigure[Transformed marginal distributions with no apparent association]
 {\includegraphics[width=.85\linewidth,clip]{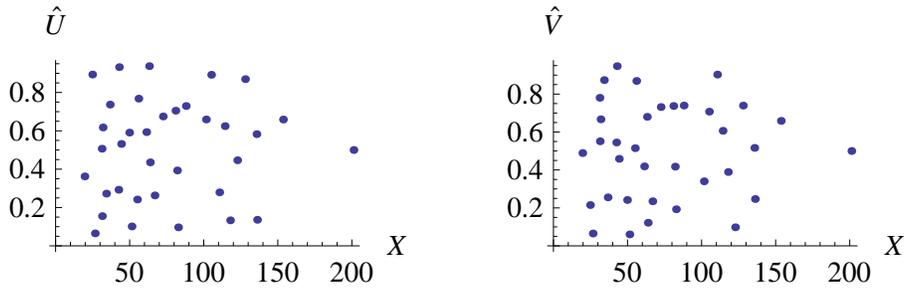}}
\caption{Illustration that the partial copula transformation of $Y_i$s into $U_i$s and $Z_i$s into $V_i$s removes marginal association. By~(\ref{eq1}), a test of independence for the $(U_i,V_i)$s is a test of conditional independence for the $(Y_i,Z_i)$s given the $X_i$s.}
 \label{scatplot1}
\end{figure}

\begin{figure}[btp]
 \centering
 {\includegraphics[width=.35\linewidth,clip]{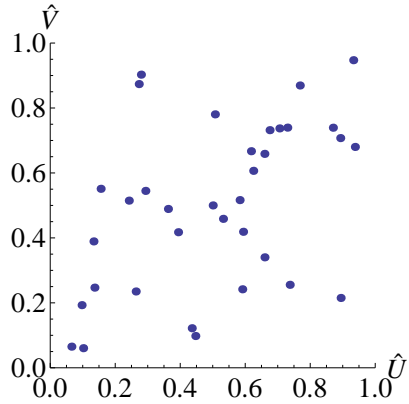}}
\caption{Partial copula transforms $(\hat U_i,\hat V_i)$ for the digoxin data. Visual inspection suggests the presence of an association for the $(\hat U_i,\hat V_i)$ which would imply a conditional assocation for the $(Y_i,Z_i)$ given the $X_i$.}
 \label{scatplot2}
\end{figure}

We now apply the aforementioned tests to a real data example, illustrating the procedure with some graphical displays.
Table~\ref{digoxin} shows data on 35 consecutive patients under
treatment for heart failure with the drug digoxin. The data are from
\citeA{halkin75}. Of medical interest is the hypothesis that digoxin
clearance is independent of urine flow controlling for creatinine
clearance, i.e., $Y\ind Z|X$.
We computed~(\ref{newobs}) using the Nadaraya-Watson estimators
\begin{eqnarray}
   \hat U_i = \hat F_{Y|X}(y|x) = \frac
   {\sum_{i=1}^n K_1(|x-X_i|/h_1)I(Y_i\le y)}
   {\sum_{i=1}^n K_1(|x-X_i|/h_1)}
   \label{kern1}
\end{eqnarray}
and
\begin{eqnarray}
   \hat V_i = \hat F_{Z|X}(z|x) = \frac
   {\sum_{i=1}^n K_2(|x-X_i|/h_2)I(Z_i\le z)}
   {\sum_{i=1}^n K_2(|x-X_i|/h_2)}
   \label{kern2}
\end{eqnarray}
where $h_t>0$ is a bandwidth and $K_t$ a kernel function ($t=1,2$).
We used a standard normal kernel and chose the bandwidth using Silverman's rule of thumb: $h_1=h_2=1.06\sigma_X n^{-1/5}=22.48$. See \cite{yj98,hwy99} and \citeA{koenker05} for other methods of estimation and bandwidth selection.
In Figure~\ref{scatplot1}(b) scatter plots are given of the pairs $(X_i,\hat U_i)$ and $(X_i,\hat V_i)$.
A visual inspection of both pictures confirms that the effect of $X$ has been removed, that is, independence seems to hold.
A scatterplot of the $(\hat U_i,\hat V_i)$ is given in
Figure~\ref{scatplot2} and some dependence is apparent, which can be tested for in various ways;
in Table~\ref{pvals}, $p$-values are given for tests based on several statistics.
The $p$-values were approximated by the permutation test using $10^5$ resamples. The simulations in Section~\ref{sec-sim} indicate that the tests are likely to be slightly liberal. Even taking that into account, there still appears to be good evidence for lack of conditional independence.

\begin{table}
\begin{center}
\begin{tabular}{rrr|rrr|rrr}
   \hline
   $X$ & $Y$ & $Z$ &   $X$ & $Y$ & $Z$ &   $X$ & $Y$ & $Z$\\
 19.5 & 17.5& 0.74 &  51.3 & 22.7& 0.33 &  101.5& 110.6& 1.38\\
 24.7 & 34.8& 0.43 & 55.0 & 30.7& 0.80 & 105.0& 114.4& 1.85\\
 26.5 & 11.4& 0.11 & 55.9 & 42.5& 1.02 & 110.5&  69.3& 2.25\\
 31.1 & 29.3& 1.48 & 61.2 & 42.4& 0.56 &  114.2&  84.8& 1.76\\
 31.3 & 13.9& 0.97 & 63.1 & 61.1& 0.93 &  117.8&  63.9& 1.60\\
 31.8 & 31.6& 1.12 & 63.7 & 38.2& 0.44 &  122.6&  76.1& 0.88\\
 34.1 & 20.7& 1.77 & 66.8 & 37.5& 0.50 & 127.9& 112.8& 1.70\\
 36.6 & 34.1& 0.70 & 72.4 & 50.1& 0.97 &  135.6& 82.2& 0.98\\
 42.4 & 25.0& 0.93 & 80.9 & 50.2& 1.02 &  136.0&  46.8& 0.94\\
 42.8 & 47.4& 2.50 & 82.0 & 50.0& 0.95 & 153.5& 137.7& 1.76\\
 44.2 & 31.8& 0.89 &  82.7 & 31.8& 0.76 &  201.1&  76.1& 0.87\\
 49.7 & 36.1& 0.52 &  87.9 & 55.4& 1.06 &\\
 \hline
\end{tabular}
\caption{Digoxin clearance data. Clearances are given in
ml/min/1.73m$^2$, urine flow in ml/min. Source: Halkin et al.\
(1975).
\newline Note: $X$= Creatinine clearance, $Y$= digoxin clearance, $Z$= urine flow. \label{digoxin}}
\end{center}
\end{table}



\begin{table}
\begin{center}
\begin{tabular}{cr}
   \hline
 Test statistic       & $p$-value \\ \hline
 $\hat\rho(\hat U,\hat V)$      & .018   \\
 $\hat\tau(\hat U,\hat V)$      & .022   \\
 $\hat\Delta(\hat U,\hat V)$    & .107   \\
 $\hat\kappa(\hat U,\hat V)$    & .041   \\
 $\hat\tau^*(\hat U,\hat V)$    & .055   \\
 \hline
\end{tabular}
\caption{Test statistics and associated $p$-values for the hypothesis $Y\ind Z|X$ for the digoxin data}\label{pvals}
\end{center}
\end{table}

\begin{figure}[btp]
 \centering
 {\includegraphics[width=.85\linewidth,clip]{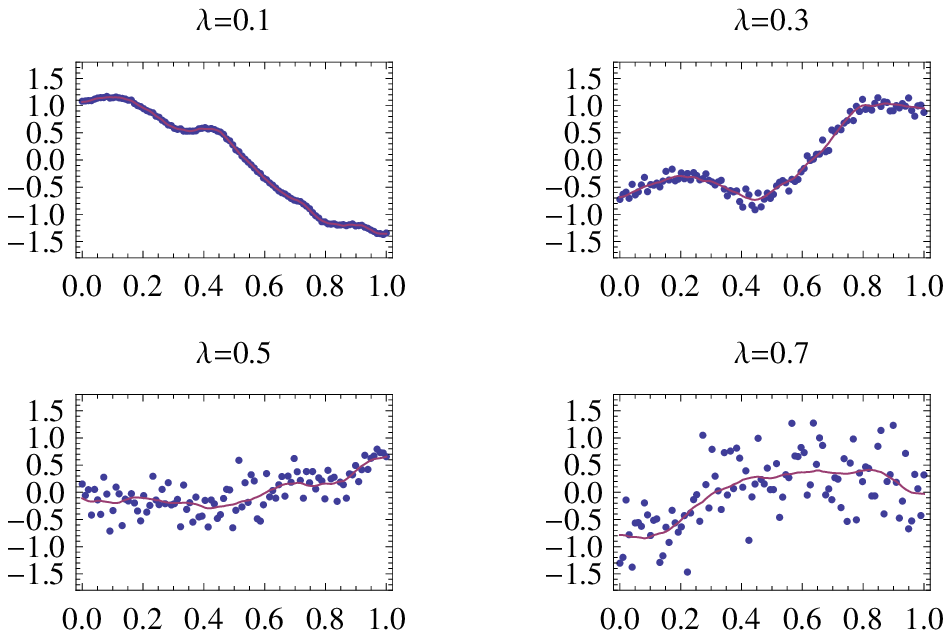}}
\caption{Randomly generated data with different noise to signal ratios $\lambda$. Here, $\sigma_0=1$ and so $\sigma_{\varepsilon}=\lambda/\sigma_0=\lambda$.}
 \label{rancurves}
\end{figure}

\section{Simulation study}\label{sec-sim}

We try to answer the following two questions regarding conditional independence tests based on the partial copula:
\begin{enumerate}
\item How good are the type I error rates?
\item How much loss of power is there compared to the corresponding unconditional test of independence?
\end{enumerate}
Data were generated according to the following model. We assumed
\begin{eqnarray*}
   Y = g(x) + \varepsilon_Y(x)  \hspace{6mm}\mbox{and}\hspace{6mm}   Z = h(x) + \varepsilon_Z(x)
\end{eqnarray*}
where the error pairs $(\varepsilon_Y(x),\varepsilon_Z(x))$ follow a bivariate normal distribution with zero mean, (partial) correlation $\rho_{YZ.X}$ and equal marginal standard deviations. The functions $g$ and $h$ were generated according to the integrated Wiener processes
\begin{eqnarray}
   g(x) = \sigma_0\int_0^x W_{1t} dt
   \hspace{6mm}\mbox{and}\hspace{6mm}
   h(x) = \sigma_0\int_0^x W_{2t} dt \label{wienmod}
\end{eqnarray}
where $W_{1t}$ and $W_{2t}$ are independent Wiener processes. Note that $g(0)=h(0)=0$, which is arbitrary but does not affect the simulations. Further note that, with probability one, $g$ and $h$ are once differentiable. A `noise to signal ratio' can be defined as $\lambda=\sigma_{\varepsilon}/\sigma_0$.  In Figure~\ref{rancurves}, randomly generated (centered) curves and data are plotted with $\sigma_0=1$ and $\lambda\in\{0.1,.0.3,0.5,0.7\}$. Finally, the values of $x$ were chosen uniformly on $[0,1]$.

We simulated data according to the above model and estimated $F_{Y|X}$ and $F_{Z|X}$ using the Nadaraya-Watson estimators~(\ref{kern1} and (\ref{kern2}).
We took a standard normal kernel and chose bandwidths according to the formula
$h_1 = h_2 = 1.75\sqrt{\lambda/n}$,
which we empirically found to yield a kernel that gives a good approximation to the hat-matrix for the posterior mode for the above model.
For each simulated data set, the $p$-value for the null hypothesis of conditional independence was obtained by performing a permutation test for independence on the $\hat U_i$ and $\hat V_i$, where we used the sample correlation coefficient as a test statistic.

\begin{figure}[tbp]
\centering
 \subfigure[$n=20$]
 {\includegraphics[width=.75\linewidth,clip]{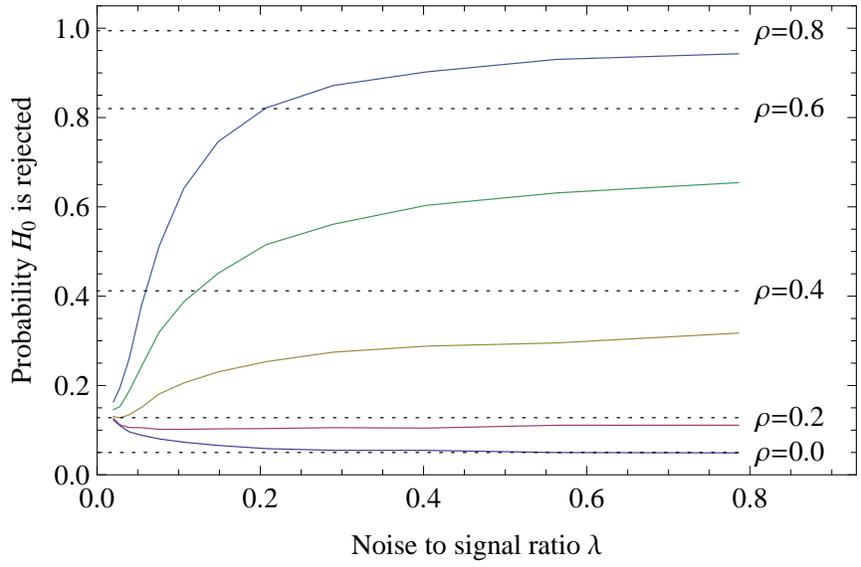}}
 \hspace{8mm}
 \subfigure[$n=100$]
 {\includegraphics[width=.75\linewidth,clip]{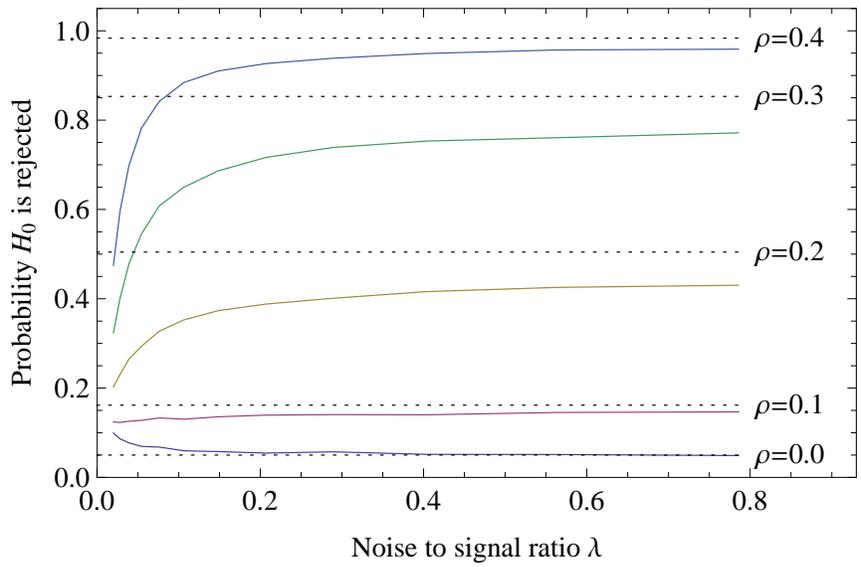}}
\caption{Power curves and Type I error rates using $\alpha=0.05$ for several values of the partial correlation $\rho$. The dotted lines give the probabilities $H_0$ is rejected for the corresponding test of marginal independence.}
 \label{powcurves}
\end{figure}

In Figure~\ref{powcurves}, it can be seen that for $n=20$ and $n=100$, unless $\lambda$ is very small, Type~I errors error rates are close to nominal and the loss of power due to conditioning seems reasonable, especially for $n=100$. For larger sample sizes the loss of power would become negligible. The method breaks down when $\lambda$ becomes very small, which is to be expected as it leads to strong overfitting.

\begin{figure}[btp]
 \centering
 \hspace*{-5mm}
 {\includegraphics[width=1.06\linewidth,clip]{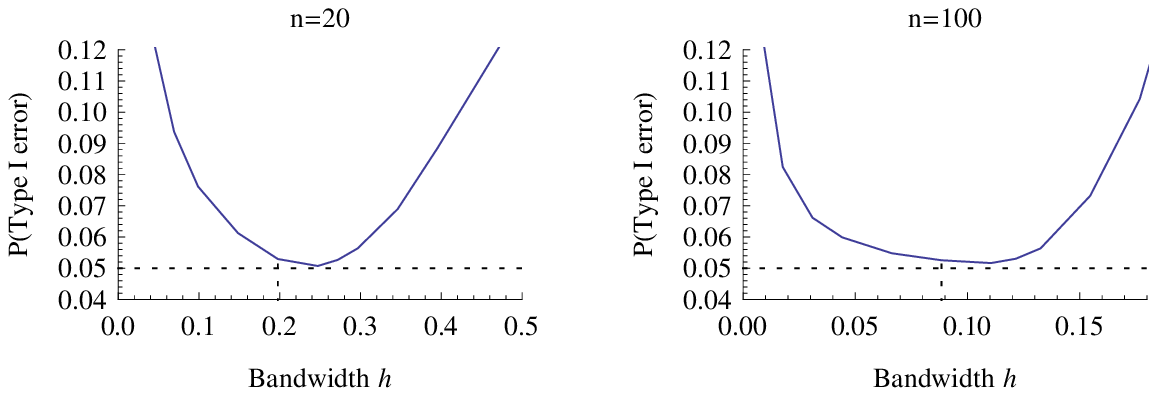}}
\caption{Type I error rates as a function of bandwidth for data generated using $\lambda=0.5$. The vertical dotted lines represent the bandwidths used in Figure~\ref{powcurves}.}
 \label{robplot}
\end{figure}

Figure~\ref{robplot} shows how the choice of bandwidth affects the Type~I error rate. In particular for $n=100$, a wide range of bandwidths give good error rates.

\section{Remarks}\label{sec-rem}

The partial copula is a tool by means of which one's favorite test of independence of appropriate form can be used for (rank) testing of conditional independence. Thus, if desired, consistency against arbitrary alternatives can be achieved by choosing an appropriate test statistic, as described in the paper.
Furthermore, our method is much simpler than competing methods.
For future research, the following questions might be investigated. Firstly, both the partial correlation and the partial copula approaches are based on removing association between response and control variables (see Section~1), but in seemingly quite different ways. This raises the question whether there is some more general approach to removing association. In particular, the case of multivariate response variables $Y$ and $Z$ remains unsolved. Secondly, in its current form the partial copula approach does not work for discrete variables. We suspect a solution to this problem can be found by considering all possible orderings of tied observations, but this needs to be investigated further. Since excellent alternative methods are available when the explanatory variable is categorical \citeA{agresti02}, this issue does not seem to be too pressing.







\bibliographystyle{apacite}

\end{document}